\documentclass[a4paper,12pt]{article}
\usepackage{amssymb}
\usepackage{amsmath}
\usepackage{amsthm}
\usepackage{indentfirst}
\usepackage{ot1patch}
\usepackage{hyperref}

\newtheorem{tw}{Theorem}[section]

\DeclareMathOperator{\dgcd}{dgcd}

\DeclareMathOperator{\cha}{char}
\DeclareMathOperator{\jac}{jac}

\DeclareMathOperator{\Irr}{Irr}
\DeclareMathOperator{\Sqf}{Sqf}
\DeclareMathOperator{\Spec}{Spec}
\DeclareMathOperator{\Prime}{Prime}
\DeclareMathOperator{\Gpr}{Gpr}
\DeclareMathOperator{\Rdl}{Rdl}

\author{Piotr J\k{e}drzejewicz, Janusz Zieli\'nski}
\title{An approach to the Jacobian Conjecture\\
in terms of irreducibility}
\date{}

\begin{document}

\maketitle

\begin{abstract}
We present some motivations and discuss various aspects
of an approach to the Jacobian Conjecture in terms of
irreducible elements and square-free elements.
\end{abstract}

\begin{table}[b]\footnotesize\hrule\vspace{1mm}
\end{table}

\section{Introduction}

The Jacobian Conjecture is one of the most important
open problems stimulating modern mathematical research
(\cite{Smale}).
Its long history is full of equivalent formulations
and wrong proofs.
In this article we give a survey of a new purely algebraic
approach to the Jacobian Conjecture in terms of irreducible
elements and square-free elements, based mainly on:
one of the authors' paper \cite{charkeller}, de Bondt and
Yan's paper \cite{BondtYan}, our paper \cite{analogs},
and our joint paper with Matysiak~\cite{factprop}.

\medskip

Let $k$ be a field of characteristic zero.
By $k[x_1,\dots,x_n]$ we denote the $k$-algebra
of polynomials in $n$ variables.
Given polynomials $f_1,\dots,f_n\in k[x_1,\dots,x_n]$,
by $\jac(f_1,\dots,f_n)$ we denote their Jacobian determinant:
$$\jac(f_1,\dots,f_n)=
\left|\begin{array}{ccc}
\frac{\partial f_1}{\partial x_1}&\cdots&
\frac{\partial f_1}{\partial x_n}\\
\vdots&&\vdots\\
\frac{\partial f_n}{\partial x_1}&\cdots&
\frac{\partial f_n}{\partial x_n}
\end{array}\right|.$$

\medskip

The Jacobian Conjecture was stated by Keller (\cite{Keller})
in 1939 for polynomials with integer coefficients.
For arbitrary field $k$ of characteristic zero it asserts
the following:
\begin{center}
{\em "If polynomials $f_1,\dots,f_n\in k[x_1,\dots,x_n]$
satisfy the Jacobian con-\\ dition
$\jac(f_1,\dots,f_n)\in k\setminus\{0\}$,
then $k[f_1,\dots,f_n]=k[x_1,\dots,x_n]$."}
\end{center}
It is known (\cite{ConnelvdDries})
that formulations of the Jacobian Conjecture for various fields
of characteristic zero (as well as for ${\mathbb Z}$)
are equivalent to each other.
The conjecture can be expressed in terms of $k$-endomorphisms of
the polynomial ring $k[x_1,\dots,x_n]$:
\begin{center}
{\em "If a $k$-endomorphism $\varphi$ of $k[x_1,\dots,x_n]$
satisfies the Jacobian con-\\ dition
$\jac(\varphi(x_1),\dots,\varphi(x_n))\in k\setminus\{0\}$,
then it is an automorphism."}
\end{center}
For more information on the Jacobian Conjecture we refer the reader
to van den Essen's book \cite{EssenBook}.

\medskip

A primary motivation of our approach can be found in a question
of van den Essen and Shpilrain from 1997 (\cite{vdEssenShpilrain},
Problem 1), whether if a $k$-endomorphism of $k[x_1,\dots,x_n]$
over a field $k$ of characteristic zero
maps variables to variables, then it is an automorphism.
A polynomial $f\in k[x_1,\dots,x_n]$ is called a {\em variable}
if there exist polynomials $f_2,\dots,f_n\in k[x_1,$ $\dots,$
$x_n]$ such that $k[f,f_2,\dots,f_n]=k[x_1,\dots,x_n]$.
A positive solution of this problem was obtained by
Jelonek (\cite{Jelonek1}, \cite{Jelonek2}).
In 2006 Bakalarski proved an analogical fact for
irreducible polynomials over~$\mathbb{C}$
(\cite{Bakalarski}, Theorem~3.7, see also \cite{Adjamagbo}).
Namely, he proved that a~complex polynomial
endomorphism is an automorphism if and only if
it maps irreducible polynomials to irreducible polynomials.
One of the authors in 2013 obtained a characterization
of $k$-endomorphisms of $k[x_1,\dots,x_n]$
satisfying the Jacobian condition as mapping
irreducible polynomials to square-free polynomials
(\cite{charkeller}, Theorem 5.1).
This fact has been further generalized by de Bondt and Yan:
they proved that mapping square-free polynomials
to square-free ones is also equivalent to the Jacobian condition
(\cite{BondtYan}, Corollary 2.2).

\medskip

We present our generalization of the Jacobian Conjecture
for $r$ polynomials $f_1,\dots,f_r\in k[x_1,\dots,x_n]$,
where $k$ is a field of characteristic zero and $r\leqslant n$:
if all jacobians (with respect to $r$ variables)
are relatively prime, then $k[f_1,\dots,f_r]$
is algebraically closed in $k[x_1,\dots,x_n]$ (\cite{analogs}).
Then we present equivalent versions of this generalized
Jacobian condition in terms of the mentioned $k$-subalgebra:
all irreducible (resp.\ square-free) elements of $k[f_1,\dots,f_r]$
are square-free in $k[x_1,\dots,x_n]$ (\cite{analogs}, Theorem 2.4).
Recall that an element $a\in R$ is called {\em square-free} if
it cannot be presented in the form $a=b^2c$, where $b,c\in R$
and $b$ is non-invertible.
It is reasonable to consider such properties in a general case,
e.g.\ for subrings of unique factorization domains.
In this case the property that square-free elements of a subring are
square-free in the whole ring can be expressed in some factorial form
(\cite{analogs}, Theorem 3.4).
At the end we discuss possible directions of future research.

\section{Freudenburg's lemma and its generalizations}

A motivation of the main preparatory fact (Theorem \ref{t24} below)
comes from generalizations of the following lemma of Freudenburg
from \cite{FreudenburgNote}.

\begin{tw}
\label{t21}
{\em (Freudenburg's Lemma)} \\
Given a polynomial $f\in \mathbb{C}[x,y]$,
let $g\in \mathbb{C}[x,y]$ be an irreducible
non-constant common factor of
$\frac{\partial f}{\partial x}$ and $\frac{\partial f}{\partial y}$.
Then there exists $c\in \mathbb{C}$ such that $g$ divides $f+c$.
\end{tw}

The assertion of the above lemma can be strengthened in a way that
if $g$ is irreducible, then
$$g\mid\frac{\partial f}{\partial x},\:
g\mid\frac{\partial f}{\partial y}\;
\Leftrightarrow\; g^2\mid f+c\;\mbox{for some}\;c\in\mathbb{C}.$$
Freudenburg needed this lemma to prove that if a polynomial
of the form $w(u,v)$, where $u$ and $v$ are variables,
belongs to the ring of constants of some
locally nilpotent derivation of the algebra
$\mathbb{C}[x_1,\dots,x_n]$,
then a variable also belongs to this ring.
Van den Essen, Nowicki and Tyc obtained
the following generalization of Freudenburg's Lemma
(\cite{vdEssenNowickiTyc}, Proposition~2.1).

\begin{tw}
\label{t22}
{\em (van den Essen, Nowicki, Tyc)} \\
Let $k$ be an algebraically closed field
of characteristic zero.
Let $Q$ be a prime ideal of the ring $k[x_1,\dots,x_n]$
and $f\in k[x_1,\dots,x_n]$.
If for each $i$ the partial derivative
$\frac{\partial f}{\partial x_i}$ belongs to $Q$,
then there exists $c\in k$ such that $f-c\in Q$.
\end{tw}

They noted (\cite{vdEssenNowickiTyc}, Remark~2.4) that
the assumption "$k$ is algebraically closed" cannot be dropped:
for $f=x^3+3x$ and $Q=(g)$, where $g=x^2+1$, in ${\mathbb R}[x]$
we have $g\mid f'$, but $g\nmid f-c$ for any $c\in {\mathbb R}$.
The idea of a generalization (in \cite{charoneel}) to arbitrary
field $k$ of characteristic zero was to consider, instead of $f-c$,
a polynomial $w(f)$, where $w(T)$ is irreducible.
In the mentioned example $w(T)=T^2+4$ since $g\mid f^2+4$.
In fact the Freudenburg's Lemma was generalized to the case when
the coefficient ring is a UFD of arbitrary characteristic.

\begin{tw}
\label{t23}
{\em (\cite{charoneel}, Theorem 3.1)} \\
Let $K$ be a unique factorization domain,
let $Q$ be a prime ideal of $K[x_1,$ $\dots,$ $x_n]$.
Consider a polynomial $f\in K[x_1,\dots,x_n]$ such that
$\frac{\partial f}{\partial x_i}\in Q$ for $i=1,\dots,n$.

\medskip

\noindent
{\bf a)}
If $\cha K=0$, then there exists an irreducible polynomial
$w(T)\in K[T]$ such that $w(f)\in Q$.

\medskip

\noindent
{\bf b)}
If $\cha K=p>0$, then there exist $b,c\in K[x_1^p,\dots,x_n^p]$
such that \mbox{$\gcd(b,c)=1$}, $b\not\in Q$ and $bf+c\in Q$.
\end{tw}

As a consequence we obtain (see \cite{charoneel}, Proposition 3.3)
that if $k$ is an arbitrary field of characteristic zero,
$f,g\in k[x_1,\dots,x_n]$ and $g$ is irreducible, then
$$g\mid\frac{\partial f}{\partial x_i}\;\mbox{for}\;
i=1,\dots,n\;\Leftrightarrow\;
g^2\mid w(f)\;\mbox{for some irreducible}\;w(T)\in k[T].$$

\medskip

A generalization of Freudenburg's Lemma to an arbitrary number
of polynomials over a field of characteristic zero
was obtained in \cite{analogs}.
Denote by $\jac^{f_1,\dots,f_r}_{x_{j_1},\dots,x_{j_r}}$
the Jacobian determinant of polynomials $f_1$, $\dots$, $f_r$
with respect to $x_{j_1}$, $\dots$, $x_{j_r}$.

\begin{tw}
\label{t24}
{\em (\cite{analogs}, Theorem 2.3)} \\
Let $k$ be a field of characteristic zero,
let $f_1,\dots,f_r\in k[x_1,$ $\dots,x_n]$ be arbitrary
polynomials, where $r\in\{1,\dots,n\}$,
and let $g\in k[x_1,$ $\dots,$ $x_n]$ be an irreducible polynomial.
The following conditions are equivalent:

\medskip

\noindent
{\rm (i)} \
$g\mid\jac^{f_1,\dots,f_r}_{x_{j_1},\dots,x_{j_r}}$ for every
$j_1,\dots,j_r\in\{1,\dots,n\}$,

\medskip

\noindent
{\rm (ii)} \
$g^2\mid w(f_1,\dots,f_r)$ for some irreducible polynomial
$w\in k[x_1,\dots,x_r]$,

\medskip

\noindent
{\rm (iii)} \
$g^2\mid w(f_1,\dots,f_r)$ for some square-free polynomial
$w\in k[x_1,\dots,x_r]$.
\end{tw}

The proof is based on the methods of proofs of earlier
special cases: Theorem 4.1 from \cite{charkeller}
and de Bondt and Yan's Theorem 2.1 from \cite{BondtYan}.

\medskip

Note also that a positive characteristic analog of Freudenburg's Lemma
for $r$ polynomials in $n$ variables was obtained in \cite{charpbases}.
It was connected with a characterization of $p$-bases
of rings of constants with respect to polynomial derivations.

\section{A characterization of Keller maps}

In this section we present the main result of \cite{charkeller}
and its substantial extension by de Bondt and Yan from
\cite{BondtYan}.
Note the following consequence of Theorem \ref{t24}
in the case $r=n$.

\begin{tw}
\label{t31}
{\em (\cite{charkeller}, Corollary 4.2, \cite{BondtYan}, Corollary 2.2)} \\
Let $k$ be a field of characteristic zero.
For arbitrary polynomials $f_1,\dots,f_n\in k[x_1,\dots,x_n]$
the following conditions are equivalent:

\medskip

\noindent
{\rm (i)} \
$\jac(f_1,\dots,f_n)\in k\setminus\{0\}$,

\medskip

\noindent
{\rm (ii)} \
for every irreducible polynomial $w\in k[x_1,\dots,x_n]$
the polynomial $w(f_1,$ $\dots,$ $f_n)$ is square-free,

\medskip

\noindent
{\rm (iii)} \
for every square-free polynomial $w\in k[x_1,\dots,x_n]$
the polynomial $w(f_1,$ $\dots,$ $f_n)$ is square-free.
\end{tw}

The above equivalence can be expressed as a characterization
of endomorphisms satisfying the Jacobian condition
analogous to the characterization of automorphisms
from Bakalarski's theorem (\cite{Bakalarski}, Theorem~3.7).

\begin{tw}
\label{t32}
{\em (\cite{charkeller}, Theorem 5.1, \cite{BondtYan}, Corollary 2.2)} \\
Let $k$ be a field of characteristic zero.
Let $\varphi$ be a $k$-endomorphism of the algebra
of polynomials $k[x_1,\dots,x_n]$.
The following conditions are equivalent:

\medskip

\noindent
{\rm (i)} \
$\jac(\varphi(x_1),\dots,\varphi(x_n))\in k\setminus\{0\}$,

\medskip

\noindent
{\rm (ii)} \
for every irreducible polynomial $w\in k[x_1,\dots,x_n]$
the polynomial $\varphi(w)$ is square-free,

\medskip

\noindent
{\rm (iii)} \
for every square-free polynomial $w\in k[x_1,\dots,x_n]$
the polynomial $\varphi(w)$ is square-free.
\end{tw}

In this way we obtain a new equivalent formulation
of the Jacobian Conjecture for an arbitrary field $k$
of characteristic zero:
\begin{quote}
{\em "Every $k$-endomorphism of $k[x_1,\dots,x_n]$
mapping square-free\\ polynomials
to square-free polynomials is an automorphism."}
\end{quote}

\medskip

There is a natural question if there exists a non-trivial example
of an endomorphism satisfying condition (ii): such that
$\varphi(w)$ is reducible for some irreducible $w$.
An affirmative answer to this question is equivalent to the negation
of the Jacobian Conjecture (\cite{charkeller}, Section 6, Remark 1).

\section{A generalization of the Jacobian Conjecture}

In \cite{analogs} we generalized the Jacobian Conjecture
in the following way
(recall that $\jac^{f_1,\dots,f_r}_{x_{j_1},\dots,x_{j_r}}$
denotes the Jacobian determinant of polynomials $f_1,\dots,f_r$
with respect to $x_{j_1}$, $\dots$, $x_{j_r}$).

\medskip

\noindent
{\bf\boldmath Conjecture $\text{JC}(r,n,k)$.}
{\em For arbitrary polynomials
$f_1,\dots,f_r\in k[x_1,$ $\dots,$ $x_n]$,
where $k$ is a field of characteristic zero
and $r\in\{1,\dots,n\}$,
if $$\gcd\big(\jac^{f_1,\dots,f_r}_{x_{j_1},\dots,x_{j_r}},\;
1\leqslant j_1<\ldots<j_r\leqslant n\big)\in k\setminus\{0\},$$
then $k[f_1,\dots,f_r]$ is algebraically closed
in $k[x_1,\dots,x_n]$.}

\medskip

Recall that by Nowicki's characterization the above assertion
means that $R$ is a ring of constants of some $k$-derivation
of $k[x_1,\dots,x_n]$ (\cite{Nrings}, Theorem~5.5,
\cite{Polder}, Theorem~4.1.5, \cite{DaigleBook}, 1.4).

\medskip

We have:

\medskip

\noindent
-- $\text{JC}(r,n,k)$ implies the ordinary Jacobian Conjecture
for $r$ polynomials in $r$ variables over $k$
(\cite{analogs}, Lemma 1.1),

\medskip

\noindent
-- $\text{JC}(1,n,k)$ is true (Ayad 2002, \cite{Ayad}, Proposition 14,
see also \cite{closed}, a remark before Proposition 4.2),

\medskip

\noindent
-- the reverse implication in $\text{JC}(r,n,k)$
need not to be true if $r<n$, we may take for example
$f_1=x_1^2x_2$, $f_2=x_3$, $\dots$, $f_r=x_{r+1}$
(\cite{analogs}, Remark 1.2).

\section{Analogs of Jacobian conditions for subrings}

In this section we present equivalent versions of the generalized
Jacobian condition from conjecture $\text{JC}(r,n,k)$ in terms of
irreducible elements as well as square-free elements.
It is useful to introduce (following \cite{jaccond})
the notion of a "differential gcd" for $r$ polynomials
$f_1,\dots,f_r\in k[x_1,\dots,x_n]$, where $r\in\{1,\dots,n\}$:
$$\dgcd(f_1,\dots,f_r)=
\gcd\big(\jac^{f_1,\dots,f_r}_{x_{j_1},\dots,x_{j_r}},\;
1\leqslant j_1<\ldots<j_r\leqslant n\big).$$

\medskip

The next theorem is a consequence of Theorem \ref{t24}
(for arbitrary $r$).

\begin{tw}
\label{t51}
Let $k$ be a field of characteristic zero.
Assume that polynomials $f_1,\dots,f_r\in k[x_1,\dots,x_n]$
are algebraically independent over $k$, where $r\in\{1,\dots,n\}$.
Then the following conditions are equivalent:

\medskip

\noindent
{\rm (i)} \
$\dgcd(f_1,\dots,f_r)\in k\setminus\{0\}$,

\medskip

\noindent
{\rm (ii)} \
for every irreducible polynomial $w\in k[x_1,\dots,x_r]$
the polynomial $w(f_1,$ $\dots,$ $f_r)$ is square-free,

\medskip

\noindent
{\rm (iii)} \
for every square-free polynomial $w\in k[x_1,\dots,x_r]$
the polynomial $w(f_1,$ $\dots,$ $f_r)$ is square-free.
\end{tw}

Note that under the assumptions of the above theorem,
a polynomial $w\in k[x_1,\dots,x_r]$ is irreducible
(square-free) if and only if $w(f_1,\dots,f_r)$
is an irreducible (square-free) element of $k[f_1,\dots,f_r]$.
This allows us to express the above conditions
in terms of the sets of irreducible elements ($\Irr$)
and square-free elements ($\Sqf$) of the respective rings.

\begin{tw}
\label{t52}
{\em (\cite{analogs}, Theorem 2.4)} \\
Let $A=k[x_1,\dots,x_n]$, where $k$ is a field
of characteristic zero.
Assume that $f_1,\dots,f_r\in A$ are algebraically
independent over $k$, where $r\in\{1,\dots,n\}$.
Put $R=k[f_1,\dots,f_r]$.
Then the following conditions are equivalent:

\medskip

\noindent
{\rm (i)} \
$\dgcd(f_1,\dots,f_r)\in k\setminus\{0\}$,

\medskip

\noindent
{\rm (ii)} \
$\Irr R\subset\Sqf A$,

\medskip

\noindent
{\rm (iii)} \
$\Sqf R\subset\Sqf A$.
\end{tw}

Therefore we may consider conditions (ii) and (iii)
in a general case, when $A$ is a domain (a commutative ring
with unity without zero divisors) and $R$ is a subring of $A$,
and we may call them analogs of the Jacobian condition~(i).
Conjecture $\text{JC}(r,n,k)$ motivated us to state the following
question (\cite{analogs}, Section 3).

\bigskip

\noindent
{\bf A general question.}
{\em Let $R$ be a subring of a domain $A$ such that
$$\Irr R\subset\Sqf A\quad\mbox{or}\quad\Sqf R\subset\Sqf A.$$
When $R$ is algebraically closed in $A$?}

\bigskip

In particular, the ordinary Jacobian Conjecture for $r=n$,
$A=k[x_1,$ $\dots,$ $x_n]$, where $\cha k=0$, asserts that
if $f_1,\dots,f_n\in A$ are algebraically independent over $k$,
$R=k[f_1,\dots,f_n]$ and $\Sqf R\subset\Sqf A$, then $R=A$.

%

\medskip

In order to understand more general context of conditions
$\Irr R\subset\Sqf A$ and $\Sqf R\subset\Sqf A$
when $R$ is a subring of a domain $A$, we can inscribe them
into the following diagram of implications
(\cite{factprop}, Proposition 3.3).

\renewcommand{\arraycolsep}{0.5mm}
\begin{small}$$\begin{array}{ccccccc}
\Irr R\subset\Irr A & \Rightarrow & \Prime R\subset\Irr A &
\Leftarrow & \Prime R\subset\Prime A & \Leftarrow &
\forall_{I\in\Spec R} AI\in\Spec A \\
\Downarrow & & \Downarrow & & \Downarrow & & \Downarrow \\
\Irr R\subset\Sqf A & \Rightarrow & \Prime R\subset\Sqf A &
\Leftarrow & \Prime R\subset\Gpr A & \Leftarrow &
\forall_{I\in\Spec R} AI\in\Rdl A \\
\Uparrow & & \Uparrow & & \Uparrow & & \Uparrow \\
\Sqf R\subset\Sqf A & \Rightarrow & \Gpr R\subset\Sqf A &
\Leftarrow & \Gpr R\subset\Gpr A & \Leftarrow &
\forall_{I\in\Rdl R} AI\in\Rdl A
\end{array}$$\end{small}

\noindent
By $\Prime R$ we have denoted the set of all prime elements
of $R$, by $\Gpr R$ the set of (single) generators of principal
radical ideals of $R$, and by $\Rdl R$ (following
\cite{BelluceDiNiolaFerraioli}, p.\ 68) the set of radical
ideals of $R$.

\section{Factorial properties}

Now we discuss factorial properties connected with inclusions
$\Irr R\subset\Irr A$ and $\Sqf R\subset\Sqf A$,
where $R$ is a subring of a unique factorization domain $A$.

\medskip

Recall that a subring $R$ of a domain $A$ such that
for every $x,y\in A$:
$$xy\in R\setminus\{0\}\;\Rightarrow\;x,y\in R,$$
is called {\em factorially closed}.
Rings of constants of locally nilpotent derivations
in domains of characteristic zero are factorially closed
(see \cite{FreudenburgBook} and \cite{DaigleBook} for details).
Note that according only to the multiplicative structure,
a submonoid of a (commutative cancelative) monoid
satisfying the above condition is called {\em divisor-closed}
(\cite{GeroldingerHalterKoch}).
Denote by $R^{\ast}$ the set of all invertible elements of a ring $R$.
It is well known that a subring $R$ of a unique factorization
domain $A$ such that $R^{\ast}=A^{\ast}$ is factorially closed
in $A$ if and only if $\Irr R\subset\Irr A$ (see \cite{analogs},
Lemma 3.2).

\medskip

Under natural assumptions we can express also the condition
$\Sqf R\subset\Sqf A$ in a form of factoriality.
If $R$ is a domain, by $R_0$ we denote its field of fractions.

\begin{tw}
\label{t61}
{\em (\cite{analogs}, Theorem 3.4)} \\
Let $A$ be a unique factorization domain.
Let $R$ be a subring of $A$ such that
$R^{\ast}=A^{\ast}$ and $R_0\cap A=R$.
The following conditions are equivalent:

\medskip

\noindent
{\rm (i)} \
$\Sqf R\subset\Sqf A$,

\medskip

\noindent
{\rm (ii)} \
for every $x\in A$, $y\in\Sqf A$,
if $x^2y\in R\setminus\{0\}$, then $x,y\in R$.
\end{tw}

If $A$ is a UFD, then a subring $R$ of $A$
that fulfills condition (ii) of Theorem \ref{t61}
we will call {\em square-factorially closed} in~$A$.
Condition (ii) has an advantage over condition (i)
since it does not involve square-free elements of $R$.
For example, one can define the square-factorial closure of
a subring $R$ in $A$ as an intersection of all square-factorially
closed subrings of $A$ containing $R$.

\medskip

There arise two questions concerning the condition
$\Irr R\subset\Sqf A$ in the case when $A$ is a UFD.
Firstly, is it equivalent to $\Sqf R\subset\Sqf A$
under some natural assumptions (like $R^{\ast}=A^{\ast}$)?
If such equivalence does not hold in general,
can the condition $\Irr R\subset\Sqf A$ be expressed
in a form of factoriality, similarly to the above theorem?

\medskip

The notion of square-factorial closedness
is relevant to thoroughly studied notion of root closedness.
Recall that a subring $R$ of a ring $A$ is called
{\em root closed} in $A$ if the following implication:
$$x^n\in R\;\Rightarrow\;x\in R$$
holds for every $x\in A$ and $n\geqslant 1$.

\begin{tw}
\label{t62}
{\em (\cite{analogs}, Theorem 3.6)} \\
Let $A$ be a unique factorization domain.
Let $R$ be a subring of $A$ such that
$R^{\ast}=A^{\ast}$ and $R_0\cap A=R$.
If $R$ is square-factorially closed in $A$,
then $R$ is root closed in $A$.
\end{tw}

An interesting task would be to investigate
whether square-factorial closedness is
stable under various operations and extensions.
Such kind of results were obtained for example
for root closedness
(see \cite{Anderson}, \cite{Angermuller},
\cite{BrewerCostaMcCrimmon}, \cite{Watkins}).
The latter is stable for instance
under homogeneous grading
and under passages to polynomial extension,
to power series extension, to rational
functions extension, to semigroup
ring $R[X;\Gamma]$, where $\Gamma$ is
torsionless grading monoid.
If for square-factorial closedness some property
would not be valid in general, then under what
additional assumptions. For example, stability
of root closure under passage to power series
extension is acquired by imposing the assumption
that a subring $R$ is von Neumann regular (see
\cite{Watkins}) or $R_0\cap A=R$ (see \cite{Angermuller})
as in Theorem \ref{t62}.
Another prospect for further research
is to obtain relationships (similarly
to Theorem \ref{t62}) of
square-factorial closedness with
other notions, such as seminormality.

\end{document}